\documentstyle[11pt]{article}
\newtheorem{lemma}{Lemma}[section]
\newtheorem{thm}[lemma]{Theorem}

\newtheorem{rmk}[lemma]{Remark}

\newtheorem{defn}[lemma]{Definition}

\font\tenmsb=msbm10
\font\sevenmsb=msbm7
\font\fivemsb=msbm5
\newfam\msbfam
\textfont\msbfam=\tenmsb
\scriptfont\msbfam=\sevenmsb
\scriptscriptfont\msbfam=\fivemsb
\def\Bbb#1{{\fam\msbfam\relax#1}}

\newcommand{\fn}[2]{\! : \! #1 \! \rightarrow \! #2}
\newcommand{\lmf}[2]{\liminf_{#1 \rightarrow #2}}
\newcommand{\lms}[2]{\limsup_{#1 \rightarrow #2}}
\newcommand{\lm}[2]{\lim_{#1 \rightarrow #2}}

\newcommand{\subsethn}{\mathrel{\mathop{\smash\subset
\vphantom{=}}\limits^\sim} }
\newcommand{\eqhn}{\mathrel{\mathop{\smash=
\vphantom{\scriptscriptstyle a}}\limits^\sim}}

\newcommand{\sigmap}{\sigma^p}

\newcommand{\wc}{\rightharpoonup}
\newcommand{\wsc}{\stackrel{*}{\rightharpoonup}}

\newcommand{\io}{\int_{\Omega}}

\newcommand{\rn}[1]{{\Bbb R}^{#1}}
\newcommand{\nat}{{\Bbb N}}
\newcommand{\hn}{{\cal H}^{N - 1}}
\newcommand{\hnc}{{\cal H}_c^{N - 1}}

\newcommand{\EE}{{\mathcal{E}}}

\begin{document}

\vspace{3cm}
\begin{center}
{\bf
QUASI-STATIC  EVOLUTION \\
IN BRITTLE FRACTURE:\\ THE CASE OF BOUNDED SOLUTIONS
}
\vspace{2cm}

Gianni {\sc Dal Maso} \\
{\scriptsize{\it S.I.S.S.A.  \\
34014 Trieste, Italia \\
email:} dalmaso@sissa.it}
\bigskip

Gilles A. {\sc Francfort} \\
{\scriptsize{\it  L.P.M.T.M., Universit\'e Paris-Nord\\
93430 Villetaneuse, France
 \\email:} francfor@galilee.univ-paris13.fr}
\bigskip

Rodica  {\sc Toader} \\
{\scriptsize{\it Dip.to di Ingegneria Civile, Universit\`a di Udine  \\
33100 Udine, Italia \\
email:} toader@dic.uniud.it}

\end{center}

\vspace{2cm}

\centerline{ \bf{Abstract} }
\medskip
\noindent
{\scriptsize The main steps of the proof of the existence result for the quasi-static evolution of cracks in brittle materials, obtained in \cite{dmft} in the vector case and for a general quasiconvex elastic energy, are presented here under the simplifying assumption that the minimizing sequences involved in the problem are uniformly bounded in~$ L^\infty$.}

\bigskip
\noindent
{\tiny{\bf Keywords:}
fracture, functions of bounded variation, geometric
measure theory, free
discontinuity problems,
{\sc  Mumford \&  Shah} functional,
quasi-static evolution.}
\bigskip

\section{Introduction}

In recent years, a variational theory of quasi-static crack growth in
a brittle solid, first proposed in \cite{fm}, has been developed on several fronts.
Its basic ingredients are few and simple.
The main
idea -- borrowed from
{\sc D. Mumford \&
J. Shah}'s approach to image segmentation \cite{ms} and close in spirit to
the original idea of {\sc A. Griffith} in  his seminal paper \cite{g} --
is
that the crack wants to quasi-statically minimize its
total energy among
all
legal competitors. In other words,
the crack $\Gamma(t)$ must minimize 
$$
\EE(t,\Gamma):={\mathcal W}(t,\Gamma) +
\hn(\Gamma)
$$
among
all $\Gamma \supset \Gamma(s),\; s<t$, where ${\mathcal W}(t,\Gamma)$
is the potential energy of the sample for the loads that are applied at time $t$, and $\hn$ is the $(N-1)$-dimensional Hausdorff measure (i.e., the length, if $N=2$, and the area, if $N=3$).

Actually, 
the evolution must be
further constrained by imposition
of a condition on the time evolution of the energy $\EE(t):=\EE(t,\Gamma(t))$, so as to recover
the propagation
criterion of {\sc A. Griffith} in the current setting. That constraint is simply
that the change in mechanical energy should exactly balance
the work of the external
loads.

Specifically, let $W$ be the elastic energy density associated to the sample
$\Omega \subset \rn{N}$, ($N= 1,\,2$, or  $3$ in the physically relevant cases).
If $u$ is any displacement field, then $W$ is either a function of $\nabla u$
(the case of finite elasticity) 
or of $\varepsilon(u)$, the symmetrized gradient of $u$ (the case of linearized 
elasticity). We assume the former and briefly comment this choice further below. 

Let $g(t)$ be the applied displacement at time $t$ on a part $\partial \Omega_d$
of the boundary  $\partial \Omega$ and let 
${\mathcal L}(t,v)$ be the work of the surface and body loads at time $t$ for
a test displacement field $v$. We do not attempt to further detail the specific kind
of surface or body loads at this point, but merely note that $\mathcal L$ {\it
cannot} be a linear map of $v$ at the current stage of the theory, as will be justified
later.  In particular, the loads cannot be displacement independent.

For a given crack $\Gamma$, a closed subset of $\overline \Omega$ with $\hn(\Gamma)<+\infty$, the displacement field
$u(t)$ is a minimizer of the potential energy
$\io W(\nabla v) dx - {\mathcal L}(t,v)$
among all kinematically admissible displacement fields, that is
all $v$'s which may be discontinuous across $\Gamma$ and satisfy the boundary condition $v = g(t)$ on $\partial \Omega_d \setminus \Gamma$.
Note that the crack may in effect debond the sample from the  applied displacement.
The resulting value of the potential energy is ${\mathcal W}(t,\Gamma)$.

Now, for a given time $t$, the crack $\Gamma(t)$ must be such that it minimizes 
the total energy $\EE(t,\Gamma)$ among all cracks that contain the prior cracks, that is
among all $\Gamma \supset \bigcup_{s<t}\Gamma(s)$.
Fracture irreversibility, that is the fact that $\Gamma(t)$ increases with $t$, is an essential feature of the
evolution process; it  is implicit in this minimality property.

Finally, the conservation of total energy must be satisfied  throughout the evolution. In the current context,
this translates into
$$
\dot \EE(t) = \int_{\partial \Omega_d\setminus \Gamma(t)}
DW(\nabla u(t))n\cdot \dot g(t)\; d\hn
- \dot {\mathcal  L}(t,u(t)).
$$

In \cite{fm}, the mechanical significance of the proposed evolution model  is investigated at length in the case
where $W$ is actually a quadratic function of $\varepsilon(u)$, the setting of linear elasticity,
$ {\mathcal  L} \equiv 0$, so that the only driving mechanism is the boundary condition $g(t)$, and
$g(t) \equiv t G$, with $G$ a fixed function. The model then palliates the major defects of 
the classical theory, most notably its inability to initiate the fracture process and to predict the crack
path as well as the time evolution of the crack along that path.

In  \cite{bc} and  \cite{bfm}
various numerical implementations of the time evolution are
attempted.
The continuous-time evolution is
replaced by a finite time
step approach. The proposed methods are shown to be both theoretically and numerically sound;
the obtained results are striking at times and certainly well beyond the boundaries of classical
fracture mechanics computations.

From the mathematical standpoint, the hurdles involved in an adequate
handling of the symmetrized gradient case -- the case where $W$ is a function of 
$\varepsilon(u)$ -- forbid at present the complete development of the theory for linearized elasticity, notably because
the ambient functional space for such a study, that is $SBD(\rn{N})$,  is only partially
understood.  The only mathematical study relevant to fracture in the context
of linearized elasticity is that of the  {\it two-dimensional setting} with a quadratic $W$, 
 and under the restrictive assumption
that  the maximal
number of connected components of the
potential cracks is  {\it a priori} known
\cite{c}.  

In the case where $\nabla u$ is considered, however, the mathematical analysis
of the proposed formulation is well under way.  The antiplane elasticity setting, that is
that where $u$ is scalar-valued, is well understood: in \cite{dmt}, existence
is shown under the same restrictions as those just detailed; then, the general quadratic case
is solved in \cite{fl}. In both settings, the  assumption  that 
$ {\mathcal  L} \equiv 0$ is essential to the analysis, because it permits to obtain $L^\infty$
estimates through the maximum principle. Since the maximum principle is not
applicable to the vector-valued case, that restriction becomes moot. In \cite{dmft}, the vector case is
analyzed under the only assumption of quasiconvexity of the energy density.  The class of 
loads $\mathcal L$ for which the minimization problem at fixed time is meaningful must ensure some 
reasonable compactness properties of the minimizing sequences, which is why a linear dependence of ${\mathcal L}(t,v)$ on $v$ is not 
admissible. See \cite{dmft} for details.

Our goal in this paper is to present the results of \cite{dmft} in a simpler situation which allows us to remove some  {\it non-trivial} technicalities  and to make the main ideas of the proof more transparent. To be definite, we
revisit the vector case under the assumption that the minimizing fields are bounded  in $L^\infty$ (uniformly in time).  We do not attempt to justify this hypothesis (which is automatically satisfied only in the scalar case) and readily agree with any potential criticism. The advantage of this assumption is that we can present our results
in a  simpler functional framework,  using the space $SBV(\rn{N};\rn{m})$ instead of 
$GSBV(\rn{N};\rn{m})$. Another advantage is that we can avoid all coerciveness assumptions on the loads, since the compactness of the minimizing sequences follows now from the $L^\infty$-bound. 

Our aim however is not to give an independent proof, but only  to describe the main arguments of \cite{dmft} in the simplest situation. For this reason we present the results only in the case of no applied loads. We still require  a few technical lemmas and refer the reader to \cite{dmft} for their proof.

Throughout, 
$SBV^p(\rn{N};\rn{m})$, $1<p<+\infty$, is the space of functions $v \in SBV(\rn{N};\rn{m})$ (see \cite{a}) such that  $\nabla v \in L^p(\rn{N};\rn{mN})$ and $\hn(S(v))<+\infty$, where $S(v)$ denotes the jump set of $v$.
We  say that a sequence 
$u^n \in SBV^p(\rn{N};\rn{m})$
$SBV^p$-converges to $u \in SBV^p(
\rn{N};\rn{m})$
 if
$$\begin{array}{c}
\nabla u^{n}\wc\nabla u\mbox{ in }
L^p(\rn{N};\rn{mN}),\\[3mm]
\hn(S(u^{n}))\mbox{ is bounded,}\\[3mm]
u^n\rightarrow u\mbox{ in } L^p(\rn{N};\rn{m}),\\[3mm]
u^{n}\wsc u\mbox{ in } L^\infty(\rn{N};\rn{m}).
\end{array}
$$
We  use the same definitions for sequences in $SBV^p(U)$,  where $U$
is an open subset of $\rn{N}$.
\begin{rmk}
\label{sbv}
If $u^n \stackrel{SBV^p}{\rightharpoonup} u$, then it is proved in \cite{a}
that, for any open set $U$,
$$\hn(S(u)\cap U)\leq  \lmf{n}{\infty}\hn(S(u^n)\cap  U).$$
Further, it is easily seen (see, e.g., \cite{dmft}, \cite{fl}) that, 
for any Borel set $E$ with  $\hn(E)<\infty$,
$$\hn(S(u)\setminus E)\leq \lmf{n}{\infty}\hn(S(u^n)\setminus E).$$
\end{rmk}

\section{Setting of the problem and statement of the results}
\label{setting}
\setcounter{equation}{0}

As mentioned in the introductory section, our analysis is restricted to the case where we act on the body only through prescribed displacements on a part $\partial \Omega_d:=\partial \Omega \setminus \partial \Omega_f$
of the boundary.

The energy density $W$ is a nonnegative quasiconvex $C^1$ function on $\rn{mN}$ which further satisfies
\begin{equation}
\label{gro}
(1/C) |\xi|^p-C\leq W(\xi)\leq C|\xi|^p+C,\quad \xi \in \rn{mN},
\end{equation}
for some constants $C\ge 1$ and $1<p<\infty$.
Note that the assumptions on $W$ immediately  imply that (see, e.g., \cite{dac})
\begin{equation}
\label{DW}
|DW(\xi)|\leq C(1+|\xi|^{p-1}),
\end{equation}
for some (possibly different) constant $C\ge 1$.

The domain $\Omega$ under consideration is assumed throughout to be bounded and Lipschitz, and the function $g$, which appears in the boundary condition
on $\partial \Omega_d$, is assumed to be defined on all of
$\rn{N}$; actually, it is taken to be in 
$W^{1,1}_{loc}([0,\infty);W^{1,p}(\rn{N};\rn{m}))$. In particular, $g$ belongs to $C^0([0,\infty);W^{1,p}(\rn{N};\rn{m}))$ and its time derivative $\dot g$ belongs to $L^1_{loc}([0,\infty);W^{1,p}(\rn{N};\rn{m}))$.

The traction-free part $\partial \Omega_f$ of the boundary $\partial \Omega$ is assumed to be closed.

We will denote throughout inclusion, up to a set of $\hn$-measure $0$, by
$\subsethn$, and set equality,  up to a set of $\hn$-measure $0$, by
$\eqhn$. A crack is a subset $\Gamma$ of $\overline{\Omega}$ with $\hn(\Gamma)<+\infty$.

 The condition that the deformation field $u$, physically defined only on $\Omega$, has a jump set contained in $\Gamma(t)$ and agrees with $g(t)$ on $\partial\Omega_d\setminus \Gamma(t)$ in the sense of traces, will be expressed in an equivalent way by defining $u$ on all of $\rn{N}$ and by requiring that $u\in SBV^p(\rn{N};\rn{m})$, $u=g(t)$ a.e.\ on $\rn{N}\setminus\overline\Omega$, and $S(u)\subsethn\Gamma(t)\cup \partial\Omega_f$. Note that, under these hypotheses $W(\nabla u)\in L^1(\Omega)$ and $DW(\nabla u) \in L^{p'}(\Omega;\rn{mN}),$ with
$p':={p/(p-1)}$.
 
 We will consider an initial crack $\Gamma_0$ and an initial deformation $u_0\in SBV^p(\rn{N};\rn{m})$ with $u_0=g(0)$ a.e.\ on $\rn{N}\setminus\overline\Omega$ and  $S(u_0)\subsethn\Gamma_0\cup \partial\Omega_f$. We assume also that the Griffith equilibrium condition is satisfied by the initial configuration $\Gamma_0$, $u_0$, i.e.,
 $u_0$ minimizes
$$
\io W(\nabla v) \; dx + \hn(S(v)\setminus(\Gamma_0\cup\partial\Omega_f))
$$
among all $v$ in $SBV^p(\rn{N};\rn{m})$ with $v = g(0) \mbox{ a.e.\ 
on }
\rn{N}\setminus \overline{\Omega}$;

In the remainder of this paper, we intend to prove the following result, under an additional uniform boundedness assumption
detailed in  {\bf (H)} below.
\begin{thm}
\label{theoreme}
There exists
a family of time dependent cracks $\Gamma(t)\subset \overline{\Omega}$, $t\ge 0$, and a
field
$u \fn{[0,\infty)\times \rn{N}}{\rn{m}} $ such that
\begin{itemize}
\item
$u(0)=u_0$ a.e.\ and $\Gamma(0)\eqhn\Gamma_0$;
\item $u(t) \in SBV^p(\rn{N};\rn{m})$  for every $t\geq 0$, so that
$DW(\nabla u(t)) \in L^{p'}(\Omega;\rn{mN})$;
\item
$ \;\Gamma(t)$ increases with $t$ and $\hn(\Gamma(t))<+\infty$ for every $t\ge 0$;
\item
$S(u(t)) \subsethn\Gamma(t)\cup\partial\Omega_f$ and $u(t) = g(t) \mbox{ a.e.\
on }
\rn{N}\setminus \overline{\Omega}$ for every $t\ge0$;
\item
for every $ t\ge0$ the deformation $u(t)$ minimizes
$$
\io W(\nabla v) \; dx + \hn(S(v)\setminus
(\Gamma(t)\cup
\partial
\Omega_f))
$$
among all $v$ in $SBV^p(\rn{N};\rn{m})$ with $v = g(t) \mbox{ a.e.\
on }
\rn{N}\setminus \overline{\Omega}$; 
\item the total energy
$$\EE(t) :=
\io W(\nabla u(t))\;dx+
\hn(\Gamma(t)\setminus \partial
\Omega_f)$$
is absolutely continuous, $ DW(\nabla u)
\cdot \nabla \dot{g} \in L^1_{loc}([0,\infty);L^1(\rn{N}))$, and
\begin{equation}
\label{expen}
\EE(t) = \EE(0)+
\int^t_0\io DW(\nabla u(s))
\cdot \nabla \dot{g}(s)\; dx\;ds
\end{equation}
for every $t> 0$.
\end{itemize}
\end{thm}

The strategy for proving the result is close to that developed in \cite{fl}
for the case of quadratic energy densities. As mentioned in the introduction, we will appeal, without repeating
the proofs, to various results in \cite{dmft}, since the present  paper favors simplicity over
completeness. 

There is no loss of generality
in restricting the study to a time interval $[0,T]$. We then choose a countable dense set
 $I_\infty$ in $[0,T]$ (containing $0$), and,
for each $n \in
\nat$,
a subset
$I_n=\{t^n_0=0 < t^n_1 < \cdots < t^n_n\}$, such that
$\{I_n\}$ form an
increasing sequence of nested sets whose union is
$I_\infty$. We set
$\Delta_n:= \sup_{k \in \{1, \ldots ,
n\}}(t^n_k-t^n_{k-1}).$ Note that
$\Delta_n \searrow
0$.

We set $\Gamma^n_0:=\Gamma_0$ and $u^n_0:=u_0$. Suppose that $u^n_j$ is defined for $j=0,1,\ldots,k-1$, and let 
$$
\Gamma^n_{k-1}:=\Gamma_0\cup 
\bigcup_{0\leq j\leq k-1}S(u^n_j).
$$
At time $t^n_k$, $k\ge1$, we define $u^n_k$ to be a minimizer for
\begin{equation}
\label{discform}
{\displaystyle \io W(\nabla v) \; dx +
\hnc(S(v)\setminus\Gamma^n_{k-1})}
\end{equation}
in $\{v \in SBV^p(\rn{N};\rn{m}): v \equiv g^n_k:=g(t^n_k) \mbox{ a.e.\ on
} \rn{N}\backslash
\overline
\Omega\}.$
In (\ref{discform}) and onward, $\hnc := \hn\lfloor \partial \Omega_f^c$, where 
$\partial \Omega_f^c:=\rn{N}\setminus \partial \Omega_f$. In other words,  $\hnc(E)=\hn(E\setminus \partial \Omega_f)$ for every Borel set $E\subset\rn{N}$.

As mentioned in the introductory section, we a priori impose that

\medskip

\begin{itemize}
\item[{\bf (H)}] {\it each problem  (\ref{discform}) has a  minimizing sequence which is bounded in $ L^\infty(\rn{N};\rn{m})$, with a bound independent  of $n$ and $k$.}
\end{itemize}
 \medskip

By a truncation argument it is easy to see that  {\bf (H)} is automatically satisfied in the scalar case $m=1$, when $g\in L^\infty([0,T] \times\rn{N})\cap
W^{1,1}([0,T];W^{1,p}(\rn{N}))$.

In view of the bound from below on $W$ and of assumption {\bf (H)}, the existence of a minimizer for
(\ref{discform}) is a straightforward iterated
application of the
$SBV$-compactness theorem (see, e.g., \cite{a}, \cite{fl}).

We then define 
$$ 
\Gamma^n(t):=\Gamma^n_k,\quad u^n(t):= u^n_k, \quad g^n(t):= g^n_k,\quad \mbox{ in }
[t^n_k,t^n_{k+1}),
$$
and  note that, for
each $t \in I_\infty$,
 $g(t) =
g_n(t)$, if $n$ is large enough.

\begin{rmk}
\label{variousmin} For every $t\in[0,T],$  $S(u^{n}(t)) \subsethn \Gamma^{n}(t)\cup \partial\Omega_f$ and $u^{n}(t)=g^{n}(t)$ a.e.\ on $\rn{N}\backslash \overline \Omega$.
Moreover,
 $u^{n}(t)$ minimizes
$$
  \io W(\nabla v)\; dx +\hnc(S(v)\setminus \Gamma^{n}(t))
$$
on $\{v \in SBV^p(\rn{N};\rn{m}): v = g^{n}(t) \mbox{ a.e.\ on }
\rn{N}\backslash \overline \Omega\}.$
In particular, with terminology from \cite{fl}, $u^{n}(t)$ is a minimizer for its own
jump set.
\end{rmk}

The construction of $\Gamma^n(t)$ and $u^n(t)$ can be viewed as a discrete time approximation of the solution to the continuous time problem. Indeed, we will also establish the following result.

\begin{thm}
\label{theorem-dc}
Consider a subsequence of $\{n\}$, independent of $t$ and still labeled $\{n\}$, such that for every $t\in[0,T]$ $\Gamma^n(t)$ $\sigmap$-converges to $\Gamma(t)$ according to  Definition~\ref{setsbv} below. 
Let $u(t)$ be a deformation field such that the pair
 $\Gamma(t), u(t)$ satisfies all the conclusions of
Theorem \ref{theoreme}.  Then for every  $t \in [0,T]$, 
$$\io W(\nabla u^n(t))\;dx \rightarrow \io  W(\nabla u(t))\;dx,$$
while
$$ \hnc(\Gamma^n(t))  \rightarrow  \hnc(\Gamma(t)).$$
\end{thm}

Note that the existence of a pair $\Gamma(t), u(t)$ is guaranteed through the proof of Theorem~\ref{theoreme}.
This latter result, a generalization in the present context of a result of 
\cite{gp}, demonstrates that the discrete time approximation provides a reasonable estimate
of both the bulk energy and the length of the crack as the discretization step becomes small.

\section{Proofs}
\label{proofq}
\setcounter{equation}{0}

As mentioned in the introduction, the proof presented here is a special case of a more general
result obtained in \cite{dmft}. In particular, we  use below the set convergenge introduced in \cite{dmft}, Section 4, under the name of $\sigma^p$-convergence,
which we now define.
\begin{defn}\label{setsbv}
We say that $\Gamma^n$ $\sigmap$-converges to $\Gamma$ if  $\Gamma^n$, $\Gamma\subset \rn{N}$, $\hn(\Gamma^n)$ is bounded uniformly with respect to $n$,  and the following conditions are satisfied:
\begin{itemize}
  \item[{\rm(a)}] if $u^j$ converges weakly to $u$ in $SBV^p(\rn{N})$ and $S(u^j)\subsethn \Gamma^{n_j}$ for some sequence $n_j\nearrow\infty$, then $S(u)\subsethn\Gamma$;
  \item[{\rm(b)}] there exist a function $u\in SBV^p(\rn{N})$ and a sequence $u^n \stackrel{SBV^p}{\rightharpoonup} u$ 
 such that $S(u)\eqhn\Gamma$ and $S(u^n)\subsethn \Gamma^n$ for every $n$.
 \end{itemize}
\end{defn}

The following compactness result  proved in Theorem 4.8 in  \cite{dmft} is central to our argument. 
Note that, in the quadratic case, one does not need to appeal to the notion of $\sigma^p$-convergence:
see \cite{fl}. Although this is true of the convex case as well, the general quasiconvex case seems
to necessitate that notion.
\begin{thm}\label{helly0}
Let $t\mapsto\Gamma^n(t)$ be a sequence of  increasing set functions defined on an interval $I$ with values contained in a bounded set  $B\subset\rn{N}$, i.e.,
$$
\Gamma^n(s)\subsethn\Gamma^n(t)\subset B\quad\hbox{for every $s,\, t\in I$ with $s<t$.}
$$
Assume that the measures
$\hn(\Gamma^n(t))$ are bounded uniformly with respect to $n$ and $t$. Then there exist a subsequence $\Gamma^{n_j}$ and an increasing set function $t\mapsto\Gamma(t)$ on $I$ such that
\begin{equation}\label{gkj}
\Gamma^{n_j}(t)\quad\sigmap\hbox{-converges to }\Gamma(t)
\end{equation}
for every $t\in I$. Furthermore, $\hn(\Gamma(t))$ is
 bounded uniformly with respect to  $t$.
\end{thm}

In all that follows, we will not relabel converging
subsequences of a given sequence, unless confusion might ensue.

\subsection{The discrete formulation}

We first derive the necessary  {\it a priori} estimates. For
some constant
$C>0$ the following holds true:
\begin{eqnarray}
&\vphantom{\displaystyle \bigcup_{\tau\leq t}}
 \|\nabla u^n(t)\|_{L^p(\rn{N};\rn{mN})}\leq C,
 \label{ape-1}
\\
&{\displaystyle \hn(\Gamma^n(t))= 
 \hn\Big(\Gamma_0\cup\bigcup_{\tau\leq t}S(u^n(\tau))\Big)}\leq C,
 \label{ape-2}
\\
&\vphantom{\displaystyle \hn\Big(}
 \|u^n(t)\|_{L^\infty(\rn{N};\rn{m})}\leq C.
 \label{ape-3}
\end{eqnarray}
Indeed, at time $t^n_k$, take $g^n_k$ as test function
for the minimality
of
$u^n_k$ in (\ref{discform}). We obtain
$$
 \io W(\nabla u^n_k)\; dx +
\hnc(S(u^n_k)\setminus
\Gamma^n_{k-1})\leq
\io W(\nabla g^n_k)\; dx,
$$
which implies, since $u^n_k \equiv g^n_k$  a.e.\ on 
$\rn{N}\backslash\overline{\Omega}$  and by virtue of the $p$-growth of $W$, that
$$
\|\nabla u^n(t)\|_{L^p(\rn{N};\rn{mN})}\leq
 C\|\nabla g^n(t)\|_{L^p(\rn{N};\rn{mN})}+C
$$
for some constant $C>0$. This proves (\ref{ape-1}), since $\|\nabla g^n(t)\|_{L^p(\rn{N};\rn{mN})}$ is bounded uniformly with respect to $t$ and~$n$.

Now, at time $t^n_{k+1}$, take $u^n_k + g^n_{k+1} -
g^n_k$ as a
test function for the minimality of $u^n_{k+1}$ in
(\ref{discform}). Since $\Phi \mapsto \io W(\Phi)\, dx$ is a 
$C^1$-map from $L^p(\rn{mN})$ into $\rn{}$ with differential
$\Psi \mapsto \io DW(\Phi) \cdot \Psi\, dx$, we
obtain, for some $\theta^n_k \in [0,1]$,
\begin{eqnarray*}
\io W(\nabla u^n_{k+1}) \;dx\!\!\!&+&\!\!\!\hnc(S(u^n_{k+1})\setminus
\Gamma^n_k )\\[3mm]
&\leq&\!\!\!
\io W(\nabla (u^n_k + g^n_{k+1} -
g^n_k)) \;dx
\\[3mm]
&=&\!\!\!
\io W(\nabla u^n_k) \;dx+
\io
DW(\nabla (u^n_k + \theta^n_k(g^n_{k+1} -
g^n_k)))  \cdot \nabla (g^n_{k+1} -
g^n_k) \; dx 
\\[3mm]
&=&\!\!\!
\io W(\nabla u^n_k)\;dx+
\int^{t^n_{k+1}}_{t^n_k} \io
 DW(\nabla u^n(s)+ \Psi^n(s)) \cdot \nabla \dot{g}(s) \; dx \;ds,
\end{eqnarray*}
where $\Psi^n \in L^\infty\big((0,T);L^p(\rn{N};\rn{mN}))$ is defined as
$$\Psi^n(s) := \theta^n_k  \int^{t^n_{k+1}}_{t^n_k}\nabla \dot{g}(\sigma) d\sigma,
\quad s \in [t^n_k,t^n_{k+1}).$$
Because  $\Delta_n\searrow 0$ and  $g(t)$ is  absolutely continuous with values in $W^{1,p}(\rn{N};\rn{m})$,
\begin{equation}
\label{remainder}
\|\Psi^n(t)\|_{L^p(\rn{N};\rn{mN})} \to 0, \; \mbox{uniformly in } t \in [0,T].
\end{equation}

Summing up the previous inequality for $k=0,\ldots,i-1$, we obtain
\begin{equation}
\label{totdisc}
\begin{array}{lll}
{\displaystyle
\io W(\nabla u^n_i)\; dx}
\!\!\!&+&\!\!\!{\displaystyle \hnc(\Gamma^n_i)}\\[5mm]\!\!\!&\leq&\!\!\!
{\displaystyle
\io W(\nabla u_0)\; dx + \hnc(\Gamma_0)}
\\[5mm]\!\!\!&&\!\!\!+ {\displaystyle \int_0^{t^n_i} \io
 DW(\nabla u^n(s)+  \Psi^n(s)) \cdot \nabla \dot{g}(s) \; dx \;ds.}
\end{array}
\end{equation}
The already established {\it a
priori} bound
(\ref{ape-1}), (\ref{remainder}), the growth estimate (\ref{DW}) on $DW$  and
H\"older's
inequality yield for every $t\in[0,T]$
\begin{equation}
\label{ouf}
 \hnc(\Gamma^n(t))\leq \io W(\nabla u(0))\;dx + \hnc(\Gamma_0) 
+C \int_0^t \|\nabla \dot{g}(s)\| _{L^p(\rn{N};\rn{mN})}ds\,,
\end{equation}
which is bounded in view of the assumed regularity of $g$. Since $\hn(\partial \Omega_f)$ is finite,
we conclude that
$$
\hn(\Gamma^n(t))\leq C.
$$
This proves (\ref{ape-2}).

The third bound (\ref{ape-3}) is an immediate consequence of  assumption {\bf (H)}.

According to Theorem \ref{helly0} applied to $\Gamma^n(t)$
and thanks to (\ref{ape-2}), there exists a subsequence of $\{n\}$,
still labeled $\{n\}$, and an increasing set function $\Gamma(t)\subset \rn{N}$
such that
\begin{equation}
\label{sig}
\Gamma^n(t) \stackrel{\sigmap}{\to} \Gamma(t)
\end{equation}
for every $t\in[0,T]$.
 Also, since  $\Gamma^n(t) \subset \overline \Omega$, we obtain that $\Gamma(t)\subsethn \overline \Omega$.
This follows from the definition of $\sigmap$-convergence and from
 Remark \ref{sbv} (applied with $U=\rn{N}\backslash\overline{\Omega}$). 

 We now set for a.e.\  $t\in[0,T]$
\begin{equation}\label{thet}
\begin{array}{rcl}\theta^{n}(t)&:=&{\displaystyle \io DW(\nabla
u^{n}(t))
\cdot \nabla
\dot{g}(t)\; dx,}
\\[4mm]\theta(t) &:=&{\displaystyle\lms{n}{\infty}\theta^{n}(t).}
\end{array}
\end{equation}

In view of the growth assumption on $DW$, the  $L^1((0,T
); L^p(\rn{N};\rn{mN}))$-regularity 
of $\nabla \dot g$, and the uniform bound (\ref{ape-1})
on  $\|\nabla u^{n}(t)\|_{L^p(\rn{N};\rn{mN})}$, {\sc Fatou}'s lemma
 immediately implies that $\theta \in L^1(0,T)$ and that
\begin{equation}
\label{theta}
\lms{n}{\infty} \int_0^t \theta^{n}(s) \;ds \leq \int_0^t\theta(s) \; ds.
\end{equation}
Furthermore, we are at liberty to extract, for a.e.\  $t\in
[0,T]$,  a 
$t$-dependent subsequence of $\theta^n$, denoted by $\theta^{n_t},$
such that
\begin{equation}
\label{thetab}
\theta(t)= \lim_{n_t\to \infty}\theta^{n_t}(t)= \lim_{n_t\to \infty}\io DW(\nabla
u^{n_t}(t))
\cdot \nabla
\dot{g}(t) \;dx.
\end{equation}

On the other hand, thanks to estimates (\ref{ape-1})--(\ref{ape-3}), we are  in a
position to apply
{\sc Ambrosio}'s
$SBV$-compactness theorem (see, e.g., \cite{a}) to
$\{u^{n_t}(t)\},$ for any $t\in[0,T]$, and to
conclude the existence of $u(t) \in SBV^p(\rn{N};\rn{m})$ 
such that, for a yet another $t$-dependent subsequence of $u^{n_t}$, 
still denoted by $u^{n_t},$
$$
\begin{array}{rcl}
\nabla u^{n_t}(t)\!\!\!&\wc&\!\!\!\nabla u(t)\mbox{ in }
L^p(\rn{N};\rn{mN}),\\[3mm]
u^{n_t}(t)\!\!\!&\rightarrow&\!\!\!u(t)\mbox{ in } L^p(\rn{N},\rn{m}),\\[3mm]
u^{n_t}(t)\!\!\!&\wsc&\!\!\!u(t)\mbox{ in } L^\infty(\rn{N};\rn{m}).
\end{array}
$$
In view of the boundedness of $\hn(S(u^{n_t}(t))$, we conclude,
following
the terminology of the introduction, that $u^{n_t}(t)$
$SBV^p$-converges to $u(t)$.

Further, (\ref{ape-1})--(\ref{ape-3}) and the lower semi-continuous character
of the $\hn$-measure with respect to
$\sigmap$-convergence (an immediate consequence of item (b) in
Definition \ref{setsbv})  imply the existence of a
constant $C$
such that,
\begin{eqnarray}
&\vphantom{\displaystyle\lmf{n_t}{\infty}}
\|\nabla u(t)\|_{L^p(\rn{N};\rn{mN})}\leq C,
\label{apet-1}
\\
&{\displaystyle \hn(\Gamma(t))}\leq\displaystyle\lmf{n_t}{\infty}\hn(\Gamma^{n_t}(t))\leq C,
\label{apet-2}
\\
&\|u(t)\|_{L^\infty(\rn{N};\rn{m})}\leq C.
\label{apet-3}
\end{eqnarray}

We now investigate the minimality properties of
$u(t)$. This is the object
of the following lemma, which is an easy
consequence of the jump transfer theorem in \cite{fl} (see Theorems  2.1, 2.8 of that reference) and of the properties of $\sigmap$-convergence.

\begin{lemma}
\label{minother}
For every $t\in[0,T]$ we have $S(u(t)) \subsethn \Gamma(t)\cup \partial\Omega_f$ and $u(t)=g(t)$ a.e.\ on $\rn{N}\backslash \overline \Omega$.
Moreover, $u(t)$ minimizes
\begin{equation}
\label{minothereq}
\io W(\nabla v) \; dx + \hnc(S(v)\setminus
\Gamma(t))
\end{equation}
on  $\{v \in SBV^p(\rn{N};\rn{m}): v = g(t)
\mbox{ a.e.\ on }
\rn{N}\backslash \overline \Omega\}$.

Further, for every $t \in [0,T]$,
\begin{equation}
\label{conven}
\io W(\nabla u^{n_t}(t))\;dx \rightarrow \io  W(\nabla u(t))\;dx,
\end{equation}
and
\begin{equation}
\label{convstr}
 DW(\nabla u^{n_t}(t))  \rightharpoonup DW(\nabla u(t)), \mbox{ weakly in } L^{p'}(\rn{N};\rn{mN}).
\end{equation}

Finally, $\theta$ defined in (\ref{thet}) lies in  $L^1(0,T)$ and
\begin{equation}
\label{thetat}
\theta(t) = \io DW(\nabla u(t)) \cdot \nabla \dot{g}(t) \; dx
\end{equation}
for almost every $t\in[0,T]$.
\end{lemma}

{\bf Proof.}
Since $u^n(t)=g^n(t)$  a.e.\ on 
$\rn{N}\backslash \overline \Omega$ we have $u(t)=g(t)$ a.e.\  on 
$\rn{N}\backslash \overline \Omega$.
That $S(u(t)) \subsethn \Gamma(t)$ is an immediate consequence of 
item (a) in Definition \ref{setsbv}, since $S(u^n(t))\subset \Gamma^n(t)$.
By item (b) in  the same definition, there exists
 $v\in SBV^p(\rn{N};\rn{m})$ with $S(v)\eqhn \Gamma(t)$ and a sequence
 $v^n\in SBV^p(\rn{N};\rn{m})$ with $S(v^n)\subsethn \Gamma^n(t)$ such that
 $v^n \stackrel{SBV^p}{\rightharpoonup} v$. We now apply the jump transfer
 theorem (Theorem  2.1 in \cite{fl}) and conclude to the existence, for an arbitrary
 element $w  \in SBV^p(\rn{N};\rn{m})$ with
  $ w = g(t) \mbox{ a.e.\ on }
\rn{N}\setminus \overline \Omega$,
 of a sequence $w^n \in SBV^p(\rn{N};\rn{m})$ such that
 \begin{equation}
 \label{tranist}
\begin{array}{c}
 w^n\equiv w = g(t) \mbox{ a.e.\ on }
\rn{N}\setminus \overline \Omega,\\[3mm]
 w^n\to  w
 \mbox{ in } L^1(\rn{N};\rn{m}),
 \\[3mm]
 \nabla w^n\to \nabla w
  \mbox{ in }L^p(\rn{N};\rn{mN}),
 \\[2.5mm]
 \hn\Big((S(w^n)\setminus S(v^n))\setminus (S(w)\setminus \Gamma(t))\Big)\longrightarrow 0.
 \end{array}
 \end{equation}
 Because $S(v^n)\subset \Gamma^n(t)$, the last inequality above {\it a fortiori} 
 implies that 
 
 \begin{equation}\label{lschn}
  \hn\Big((S(w^n)\setminus \Gamma^n(t))\setminus (S(w)\setminus \Gamma(t))\Big)
  \longrightarrow 0.
  \end{equation}

 Now, in view of Remark \ref{variousmin},
 $u^{n_t}(t)$ minimizes
$$
\io W(\nabla v) \; dx + \hnc(S(v)\setminus \Gamma^{n_t}(t))
$$
on $\{v \in SBV^p(\rn{N};\rn{m}): v = g^{n_t}(t) \mbox{ a.e.\ on }
\rn{N}\backslash \overline \Omega\}$, so that
\begin{equation}
\label{finitemin}
 \begin{array}{l}
 \displaystyle \io W(\nabla u^{n_t}(t)) \; dx \\
  \displaystyle \qquad \leq \io W(\nabla w^{n_t}(t)+ \nabla g^{n_t}(t)-\nabla g(t)) \; dx + \hnc(S(w^{n_t}(t)) \setminus \Gamma^{n_t}(t)).
 \end{array}
\end{equation}

Since $W$ is quasiconvex with $p$-growth, $u^{n_t}(t)$
$SBV^p$-converges to $u(t)$, and the sequence $\hn(S(u^{n_t}(t))$ is uniformly bounded,  Theorem 5.29 in \cite{afp} implies that
\begin{equation}
\label{helpest}
\io  W(\nabla u(t)) \; dx  \leq \lmf{n_t}{\infty}\io W(\nabla u^{n_t}(t)) \; dx.
\end{equation}
 Recalling (\ref{tranist})--(\ref{helpest}),  and the fact that 
 $\nabla g^{n_t}(t) 
 \stackrel{L^{p}(\rn{N};\rn{mN})}{\longrightarrow} \nabla g(t)$, we conclude
 that
$$
\io W(\nabla u(t)) \; dx \leq \io W(\nabla
w)\; dx +
\hnc (S(w)\setminus \Gamma(t))
$$
and obtain the
minimality result.

To prove (\ref{conven}), we apply the jump transfer theorem once again, this time to 
$u(t)$, thus obtaining a sequence $w_{n_t} \in SBV^p(\rn{N};\rn{m})$   with $ w_{n_t} \equiv g^{n_t}(t)$
a.e.\ on $\rn{N}\backslash \overline \Omega$ and such that
$$
\begin{array}{c}
\nabla w_{n_t} \rightarrow \nabla u(t) \mbox{ in }
L^p(\rn{N};\rn{mN}),\\[3mm]
\hnc(S(w_{n_t})\setminus S(u^{n_t}(t)))
\rightarrow 0.
\end{array}
$$
Since $u^{n_t}(t)$ is in particular a minimizer for its own jump set,
$$
\io W(\nabla u^{n_t}(t)) \; dx \leq \io W(\nabla
w_{n_t})\; dx +
\hnc (S(w_{n_t})\setminus S(u^{n_t}(t))).
$$
Thus
$$
\lms{n_t}{\infty} \io W(\nabla u^{n_t}(t)) \; dx \leq \io W(\nabla u(t)) \; dx,
$$
which, together with (\ref{helpest}), yields the desired result.

To prove (\ref{convstr}), we appeal to Lemma 4.11 in \cite{dmft}, which
 states in essence  that $SBV^p$-convergence of $u^{n_t}(t)$ to $u(t)$, 
together with convergence (\ref{conven}) of the 
energy, implies weak convergence of the stresses.

Finally, (\ref{thetat}) is an immediate consequence of (\ref{thetab}) and (\ref{convstr}). 
\hfill $\Box$
\bigskip

We now derive an elementary estimate on the total energy
at time $t$, that
is
\begin{equation}
\label{totalenergy}
\EE(t) :=
\io W(\nabla u(t)) \; dx+
\hnc(\Gamma(t)).
\end{equation}
We also define the corresponding total energy for
$u^n(t)$, namely
\begin{equation}
\label{totalenergyn}
\EE^n(t) :=
\io W(\nabla u^n(t)) \; dx +
\hnc(\Gamma^n(t)).
\end{equation}

The following lemma then holds.
\begin{lemma}
\label{enest}
For any $t \in [0,T],$
\begin{equation}
\label{enleq}
\EE(t) \leq \EE(0)+
 \int^t_0\io DW(\nabla u(s))
\cdot \nabla \dot{g}(s)\; dx\;ds. 
\end{equation}
\end{lemma}
{\bf Proof.}
We recall (\ref{totdisc}), namely
$$
\begin{array}{lll}
{\displaystyle
\EE^{n_t}(t)}& = &{\displaystyle \io W(\nabla u^{n_t}(t))\;dx
+ \hnc(\Gamma^{n_t}(t))}\\[4mm]&\leq&
{\displaystyle
\io W(\nabla u^{n_t}(0))\;dx + \hnc(\Gamma_0)}
\\&&\quad{\displaystyle+ \int_0^t \io
 DW(\nabla u^{n_t}(s)+  \Psi^{n_t}(s)) \cdot \nabla \dot{g}(s)\; dx \;ds.}
\end{array}
$$
and pass to the limit in ${n_t}$. The  $SBV^p$-convergence of $u^{n_t}(t)$ to $u(t)$, together
with Theorem  5.29 in \cite{afp} 
 permits us to pass to the lim-inf in the
first term of the left side of the above inequality while we appeal to (\ref{apet-2})
for  the  surface term.

Finally, in view of (\ref{remainder}), a simple argument based on the uniform continuity of 
$DW$ on compact sets, together with the already established uniform bound on
$\nabla u^n(t)$ in $L^p(\rn{N};\rn{mN})$ (cf. (\ref{ape-1})) permits to drop
$\Psi^{n_t}(s)$ in the remaining term; see Lemma~4.9 in \cite{dmft} for a complete proof in 
a more general setting. Specifically,
 for a.e. $s \in [0,t]$, Lemma~4.9 in \cite{dmft} yields
$$\io  DW(\nabla u^{n_t}(s)+  \Psi^{n_t}(s))\cdot
\nabla \dot{g}(s) \; dx-   \io DW(\nabla u^{n_t}(s))\cdot
\nabla \dot{g}(s) \; dx\stackrel{n\to \infty}{\longrightarrow} 0.$$
The growth property of $DW$, together with the uniform 
$L^p(\rn{N};\rn{mN})$-bound on $\nabla u^{n_t}$ and (\ref{remainder}), imply that
$$\int_0^t \Big(\io  DW(\nabla u^{n_t}(s)+  \Psi^{n_t}(s))\cdot
\nabla \dot{g}(s) \; dx-   \io DW(\nabla u^{n_t}(s))\cdot
\nabla \dot{g}(s) \; dx\Big)\;ds \stackrel{n\to \infty}{\longrightarrow} 0,$$
We obtain
$$
\EE(t)\leq \EE(0)+
\lmf{n_t}{\infty}\int_0^t \Big(\io DW(\nabla
u^{n_t}(s))
\cdot \nabla
\dot{g}(s)\; dx\Big)\;ds.
$$
In view of (\ref{theta}) and (\ref{thetat}), the last term in the above inequality is 
bounded from above by the announced expression.
\hfill$\Box$

\bigskip

It now remains to establish that inequality (\ref{enleq})
is actually an equality. This is the object of the following lemma.

\begin{lemma}
\label{enest2}We have
\begin{equation}
\label{engeq}
\EE(t) \geq \EE(0)+ \int^t_0\io DW(\nabla u(s))
\cdot \nabla \dot{g}(s)\; dx\;ds. 
\end{equation}
\end{lemma}
{\bf Proof.}  
We
take $v \equiv u(t) + g(s) - g(t)$ as a competitor for
$u(s)$ in the minimum problem for 
(\ref{minothereq}), and get, since $S(u(t)) \subsethn \Gamma(t)$,
\[\begin{array}{lcl}
\EE(s)&\leq& {\displaystyle\io  W(\nabla v)\;dx
+\hnc(S(u(t))\setminus \Gamma(s))
+\hnc(\Gamma(s))}
\\[8mm]&\leq&{\displaystyle
 \io W(\nabla v)\;dx) +\hnc(\Gamma(t)),}
\end{array}\]
so that, for some $\rho(s,t) \in [0,1]$,
\begin{eqnarray}
 \label{sten}
\nonumber\EE(t)-\EE(s) &\geq &\io (W(\nabla u(t))-
W(\nabla v)) \;dx \\&=&
 \io \Big[ DW\Big(\nabla u(t)+ \rho(s,t) \int_s^t \nabla \dot{g}(\tau) \; d\tau\Big) \cdot
  \int_s^t \nabla \dot{g}(\tau) \; d\tau\Big] \;dx.
\end{eqnarray}
 Consider a partition $0:=s^n_0\leq s^n_1\leq \cdots \leq s^n_{k(n)}=t$
such that
\begin{equation}\label{dn}
 \lim_{n\to \infty}\max_{1\leq i \leq k(n)}(s^n_i-s^n_{i-1})=0\,,
 \end{equation}
define 
 $$
 u_n(s):= u(s^n_{i+1})\qquad\hbox{and}\qquad
 X_n(s) := \rho(s^n_i,s^n_{i+1}) \int_{s^n_i}^{s^n_{i+1}} \nabla \dot{g}(\tau) \; d\tau\,,
 $$
for $s \in (s^n_i,s^n_{i+1}]$, and note that, 
since
$g \in W^{1,1}((0,t);W^{1,p}(\rn{N};\rn{m}))$,
\begin{equation}
\label{incratio}
\|X_n(s)\|_{L^p(\rn{N};\rn{mN})} \to 0, \mbox{ uniformly on } [0,t].
  \end{equation}
  
We apply  (\ref{sten}) for $s=s^n_i$ and $t=s^n_{i+1}$, and sum the result
for $i=0, \ldots , k(n)-1$; we obtain
$$
\EE(t)- \EE(0) \geq  \int_{0}^{t} \io DW(\nabla u_n(s) + X_n(s))\cdot
\nabla \dot{g}(s) \; dx \;ds.
$$
Recalling (\ref{incratio}), we immediately infer, using an argument similar to that
which allowed to drop the term in $\Psi^n$ in (\ref{totdisc}) for the proof
of Lemma \ref{enest} (see once again Lemma 4.9 in \cite{dmft}),
that, for a.e. $s \in [0,t]$,
$$
\io DW(\nabla u_n(s) + X_n(s))\cdot
\nabla \dot{g}(s) \; dx-   \io DW(\nabla u_n(s))\cdot
\nabla \dot{g}(s) \; dx\stackrel{n\to \infty}{\longrightarrow} 0.
$$
The growth property of $DW$, together with the uniform 
$L^p(\rn{N};\rn{mN})$-bound on $\nabla u_n$ and (\ref{incratio}), imply that
$$
\int_0^t \Big(\io DW(\nabla u_n(s) + X_n(s))\cdot
\nabla \dot{g}(s) \; dx-   \io DW(\nabla u_n(s))\cdot
\nabla \dot{g}(s) \; dx\Big)\;ds \stackrel{n\to \infty}{\longrightarrow} 0,
$$
so that
\begin{equation}
\label{bfb}
\EE(t)- \EE(0) \geq  \lms{n}{\infty}\int_{0}^{t} \io DW(\nabla u_n(s))\cdot
\nabla \dot{g}(s) \; dx \;ds.
\end{equation}

To complete the proof, we need to appeal to the following result in measure theory
(see Lemma 4.12 in \cite{dmft}).

\begin{lemma}
Let $X$ be a Banach space and $f\in L^1((0,t);X)$. Then, there exists
a sequence of subdivisions $0=s^n_0\leq s^n_1\leq \cdots \leq s^n_{k(n)}=t$, satisfying (\ref{dn}), such that
$$
\lim_{n\to\infty}\,\sum_{i=1}^{k(n)} \Big \| (s^n_i-s^n_{i-1}) f(s^n_i) -
\int_{s^n_{i-1}}^{s^n_i} f(t)\,dt \,\Big\|_X= 0.
$$
\end{lemma}

We apply this lemma to
$$
f:= (\nabla \dot{g}, \theta) \in L^1((0,t);L^p(\Omega;\rn{mN})\times \rn{}),
$$
which allows to find a sequence of subdivisions $0=s^n_0\leq s^n_1\leq \cdots \leq s^n_{k(n)}=t$, so that first $\nabla \dot{g}(s)$ is replaced by 
$$G_n(s):= \nabla \dot{g}(s^n_i),
\; s^n_{i-1}<s\leq s^n_i$$
in (\ref{bfb}), and also so that
$$
\int_0^t\io DW(\nabla u_n(s))\cdot G_n(s) \; dx\,ds \longrightarrow \int_0^t\theta(s)\; ds.
$$
In view of the expression (\ref{thetat}) for $\theta(s)$, we get the desired result.
\hfill $\Box$

\medskip

The proof of Theorem \ref{theoreme} is complete.

\bigskip

\noindent{\bf Proof of Theorem \ref{theorem-dc}.} 
First of all we observe that if the pair $\Gamma(t)$, $u(t)$ satisfies all the conclusions of Theorem~\ref{theoreme}, then the value of the integral 
\begin{equation}\label{100}
\io W(\nabla u(t))\;dx
\end{equation}
does not depend on the choice of $u(t)$. Indeed, (\ref{100}) minimizes $\io W(\nabla v)\;dx$ among all $v\in SBV^p(\rn{N};\rn{m})$ with $v=g(t)$ a.e.\ on $\rn{N}\setminus\overline \Omega$ and $S(v)\subsethn \Gamma(t)\cup\partial\Omega_f$.

Recalling the definition (\ref{totalenergyn}) of the total energy $\EE^n(t)$ and using
an argument identical to that used in the proof of Lemma \ref{enest}, we obtain,
with the function $u(t)$ constructed in the proof of  Theorem \ref{theoreme},
\begin{equation}
\label{dc4}
\begin{array}{c}
\displaystyle\lms{n}{\infty} \EE^n(t)-\EE(0) \leq \lms{n}{\infty}\int_0^t \theta^n(s)\;ds\leq \int_0^t
\theta(s)\; ds\\
\displaystyle= \int_0^t\io DW(\nabla u(s)) \cdot \nabla \dot{g}(s) \; dx\;ds,
\end{array}
\end{equation}
where we have appealed to (\ref{theta}) in deriving the last inequality and
invoked the expression (\ref{thetat}) for $\theta(t)$. Hence,   {\it a fortiori}, 
\begin{eqnarray*}
\lms{n}{\infty}\io W(\nabla u^n(t))\;dx +
\lmf{n}{\infty} \hnc(\Gamma^n(t))- \EE(0)\\
\leq \int_0^t\io DW(\nabla u(s)) \cdot \nabla \dot{g}(s) \; dx\;ds=\EE(t)-\EE(0),
\end{eqnarray*}
where the last equality follows from (\ref{expen}). In view of (\ref{apet-2})   we conclude that
\begin{equation}
\label{dc3}
\lms{n}{\infty}\io W(\nabla u^n(t))\;dx  \leq \io W(\nabla u(t)) \; dx.
\end{equation}

Now, consider a $t$-dependent sequence
$\{n_t\}$ such that
\begin{equation}
\label{dc1}
\lmf{n}{\infty} \io W(\nabla u^n(t))\;dx = \lm{n_t}{\infty} \io W(\nabla u^{n_t}(t))\; dx.
\end{equation}
The sequence $\{ u^{n_t}(t)\}$ may be assumed to $SBV^p$-converge to
some $\overline u(t) \in SBV^p(\rn{N};\rn{m})$, and, as in the proof of Lemma
\ref{minother}, $\sigmap$-convergence, together with the jump transfer theorem 
imply that, just like $u(t)$, $\overline u(t)$ minimizes
$$
\io W(\nabla v) \; dx + \hnc(S(v)\setminus
\Gamma(t))
$$
 among all $v$ in $\{v \in SBV^p(\rn{N};\rn{m}): v = g(t)
\mbox{ a.e.\ on }
\rn{N}\backslash \overline \Omega\}$, and that
$S(\overline u(t))\subsethn \Gamma(t)$.
 Thus,
 $$\io W(\nabla \overline u(t)) \; dx = \io W(\nabla u(t)) \; dx.$$
 But, by sequential weak lower semi-continuity,
$$
\io W(\nabla \overline u(t)) \; dx \leq \lmf{n}{\infty} \io W(\nabla u^n(t))\;dx,$$
 hence
\begin{equation}
\label{dc6}
\io W(\nabla u(t)) \; dx\leq \lmf{n}{\infty} \io W(\nabla u^n(t))\;dx,
\end{equation}
 which, together with (\ref{dc3}) yields 
 \begin{equation}
\label{dc7}
  \lm{n}{\infty} \io W(\nabla u^n(t))\;dx =  \io W(\nabla u(t)) \; dx.
 \end{equation}
 
 Finally, recalling (\ref{dc4}), (\ref{expen}), (\ref{apet-2}) and (\ref{dc6}), we have
 \begin{eqnarray*}\io W(\nabla u(t)) \; dx \!\!\!&+&\!\!\!  \hnc(\Gamma(t))\leq \lmf{n}{\infty} \EE^n(t) \\\!\!\!&
 \leq&\!\!\!  \lms{n}{\infty} \EE^n(t) \leq \io W(\nabla u(t)) \; dx + \!\hnc(\Gamma(t)),
 \end{eqnarray*}
 so that, 
 $$ \lm{n}{\infty} \EE^n(t)= \io W(\nabla u(t)) \; dx + \hnc(\Gamma(t))$$
 and, in view of (\ref{dc7}),
 $$\hnc(\Gamma(t))= \lm{n}{\infty} \hnc(\Gamma^n(t)).$$
 
 The proof of Theorem \ref{theorem-dc} is complete.\hfill $\Box$

\end{document}